
%
%
%
\magnification=\magstep1
\def\autori{Graziano Crasta}
\def\titolo{An existence result for non-coercive non-convex problems}
\newif\iftitlepage
\newif\ifproofmode
\proofmodefalse
\titlepagetrue
%
%
%
%
\font\twelverm=cmr12
\font\twelvei=cmmi12
\font\twelvesy=cmsy10
\font\twelvebf=cmbx12
\font\twelvett=cmtt12
\font\twelveit=cmti12
\font\twelvesl=cmsl12
\font\ninerm=cmr9
\font\ninei=cmmi9
\font\ninesy=cmsy9
\font\ninebf=cmbx9
\font\ninett=cmtt9
\font\nineit=cmti9
\font\ninesl=cmsl9
\font\ninesc = cmcsc10 at 9pt

\font\eightrm=cmr8
\font\eighti=cmmi8
\font\eightsy=cmsy8
\font\eightbf=cmbx8
\font\eighttt=cmtt8
\font\eightit=cmti8
\font\eightsl=cmsl8
\font\sixrm=cmr6
\font\sixi=cmmi6
\font\sixsy=cmsy6
\font\sixbf=cmbx6
\catcode`@=11 
\newskip\ttglue
%
\def\twelvepoint{\def\rm{\fam0\twelverm}
\textfont0=\twelverm  \scriptfont0=\ninerm
\scriptscriptfont0=\sevenrm
\textfont1=\twelvei  \scriptfont1=\ninei  \scriptscriptfont1=\seveni
\textfont2=\twelvesy  \scriptfont2=\ninesy
\scriptscriptfont2=\sevensy
\textfont3=\tenex  \scriptfont3=\tenex  \scriptscriptfont3=\tenex
\textfont\itfam=\twelveit  \def\it{\fam\itfam\twelveit}%
\textfont\slfam=\twelvesl  \def\sl{\fam\slfam\twelvesl}%
\textfont\ttfam=\twelvett  \def\tt{\fam\ttfam\twelvett}%
\textfont\bffam=\twelvebf  \scriptfont\bffam=\ninebf
\scriptscriptfont\bffam=\sevenbf  \def\bf{\fam\bffam\twelvebf}%
\tt  \ttglue=.5em plus.25em minus.15em
\normalbaselineskip=15pt
\setbox\strutbox=\hbox{\vrule height10pt depth5pt width0pt}%
\let\sc=\tenrm  \let\big=\twelvebig  \normalbaselines\rm}
%
\def\tenpoint{\def\rm{\fam0\tenrm}
\textfont0=\tenrm  \scriptfont0=\sevenrm  \scriptscriptfont0=\fiverm
\textfont1=\teni  \scriptfont1=\seveni  \scriptscriptfont1=\fivei
\textfont2=\tensy  \scriptfont2=\sevensy  \scriptscriptfont2=\fivesy
\textfont3=\tenex  \scriptfont3=\tenex  \scriptscriptfont3=\tenex
\textfont\itfam=\tenit  \def\it{\fam\itfam\tenit}%
\textfont\slfam=\tensl  \def\sl{\fam\slfam\tensl}%
\textfont\ttfam=\tentt  \def\tt{\fam\ttfam\tentt}%
\textfont\bffam=\tenbf  \scriptfont\bffam=\sevenbf
\scriptscriptfont\bffam=\fivebf  \def\bf{\fam\bffam\tenbf}%
\tt  \ttglue=.5em plus.25em minus.15em
\normalbaselineskip=12pt
\setbox\strutbox=\hbox{\vrule height8.5pt depth3.5pt width0pt}%
\let\sc=\eightrm  \let\big=\tenbig  \normalbaselines\rm}
%
\def\ninepoint{\def\rm{\fam0\ninerm}
\textfont0=\ninerm  \scriptfont0=\sixrm  \scriptscriptfont0=\fiverm
\textfont1=\ninei  \scriptfont1=\sixi  \scriptscriptfont1=\fivei
\textfont2=\ninesy  \scriptfont2=\sixsy  \scriptscriptfont2=\fivesy
\textfont3=\tenex  \scriptfont3=\tenex  \scriptscriptfont3=\tenex
\textfont\itfam=\nineit  \def\it{\fam\itfam\nineit}%
\textfont\slfam=\ninesl  \def\sl{\fam\slfam\ninesl}%
\textfont\ttfam=\ninett  \def\tt{\fam\ttfam\ninett}%
\textfont\bffam=\ninebf  \scriptfont\bffam=\sixbf
\scriptscriptfont\bffam=\fivebf  \def\bf{\fam\bffam\ninebf}%
\tt  \ttglue=.5em plus.25em minus.15em
\normalbaselineskip=11pt
\setbox\strutbox=\hbox{\vrule height8pt depth3pt width0pt}%
\let\sc=\sevenrm  \let\big=\ninebig  \normalbaselines\rm}
%
\def\eightpoint{\def\rm{\fam0\eightrm}
\textfont0=\eightrm  \scriptfont0=\sixrm  \scriptscriptfont0=\fiverm
\textfont1=\eighti  \scriptfont1=\sixi  \scriptscriptfont1=\fivei
\textfont2=\eightsy  \scriptfont2=\sixsy  \scriptscriptfont2=\fivesy
\textfont3=\tenex  \scriptfont3=\tenex  \scriptscriptfont3=\tenex
\textfont\itfam=\eightit  \def\it{\fam\itfam\eightit}%
\textfont\slfam=\eightsl  \def\sl{\fam\slfam\eightsl}%
\textfont\ttfam=\eighttt  \def\tt{\fam\ttfam\eighttt}%
\textfont\bffam=\eightbf  \scriptfont\bffam=\sixbf
\scriptscriptfont\bffam=\fivebf  \def\bf{\fam\bffam\eightbf}%
\tt  \ttglue=.5em plus.25em minus.15em
\normalbaselineskip=9pt
\setbox\strutbox=\hbox{\vrule height7pt depth2pt width0pt}%
\let\sc=\sixrm  \let\big=\eightbig  \normalbaselines\rm}
%
\def\twelvebig#1{{\hbox{$\textfont0=\twelverm\textfont2=\twelvesy
	\left#1\vbox to10pt{}\right.\n@space$}}}
\def\tenbig#1{{\hbox{$\left#1\vbox to8.5pt{}\right.\n@space$}}}
\def\ninebig#1{{\hbox{$\textfont0=\tenrm\textfont2=\tensy
	\left#1\vbox to7.25pt{}\right.\n@space$}}}
\def\eightbig#1{{\hbox{$\textfont0=\ninerm\textfont2=\ninesy
	\left#1\vbox to6.5pt{}\right.\n@space$}}}

\font\medbf=cmbx10 scaled\magstep2
\def\today{\ifcase\month\or
January\or February\or March\or April\or May\or June\or
July\or August\or September\or October\or November\or December\fi
\space\number\day, \number\year}
%
%
%
%
\nopagenumbers

\def\sqr#1#2{\vbox{
   \hrule height .#2pt
   \hbox{\vrule width .#2pt height #1pt \kern #1pt
      \vrule width .#2pt}
   \hrule height .#2pt}}

%
%
%
\newcount\vol \newcount\pag
\def\bibart#1#2#3#4#5#6#7#8{\global\vol=#5 \global\pag=#7
{\item{[{\bf\bib{#1}}]}{\ninesc #2}: {#3},
{\sl #4},%
{\ifnum\vol=0\else{{\bf\ #5},}\fi}
{#6}%
{\ifnum\pag=0\else{, p.~{#7}--{#8}}\fi}%
. \medskip}}
\def\bibprep#1#2#3#4#5{
{\item{[{\bf\bib{#1}}]}{\ninesc #2}: {#3},
{\sl #4}, {#5}. \medskip}}
\def\biblib#1#2#3#4#5#6{
{\item{[{\bf\bib{#1}}]}{\ninesc #2}: {\sl #3}, {#4}, {#5}, {#6}.
\medskip}}
\mathchardef\emptyset="001F
%
%
%
\font\sixrm=cmr6
\newcount\tagno \tagno=0		        
\newcount\thmno	\thmno=0	         	
\newcount\bibno	\bibno=0			
\newcount\chapno\chapno=0                       
\newcount\verno            
\newif\ifwanted
\wantedfalse
\newif\ifindexed
\indexedfalse
\def\ifundefined#1{\expandafter\ifx\csname+#1\endcsname\relax}
\def\Wanted#1{\ifundefined{#1} \wantedtrue \immediate\write0{Wanted #1
\the\chapno.\the\thmno}\fi}
\def\Increase#1{{\global\advance#1 by 1}}
\def\Assign#1#2{\immediate
\write1{\noexpand\expandafter\noexpand\def
 \noexpand\csname+#1\endcsname{#2}}\relax
 \global\expandafter\edef\csname+#1\endcsname{#2}}
\def\pAssign#1#2{\write1{\noexpand\expandafter\noexpand\def
 \noexpand\csname+#1\endcsname{#2}}}
\def\lPut#1{\ifproofmode\llap{\hbox{\sixrm #1\ \ \ }}\fi}
\def\rPut#1{\ifproofmode$^{\hbox{\sixrm #1}}$\fi}
%
\def\chp#1{\global\tagno=0\global\thmno=0\Increase\chapno
\Assign{#1}
{\the\chapno}{\lPut{#1}\the\chapno}}
\def\newpar#1#2{\bigskip {\noindent\bf \chp{#1}. #2} \medskip}  %
\def\thm#1{\Increase\thmno
\Assign{#1}{\the\chapno.\the\thmno}\the\chapno.\the\thmno\rPut{#1}}
\def\frm#1{\Increase\tagno
  \Assign{#1}{\the\chapno.\the\tagno}\lPut{#1}{\the\chapno.\the\tagno}}

\def\bib#1{\Increase\bibno
\Assign{#1}{\the\bibno}\lPut{#1}{\the\bibno}}
\def\pgp#1{\pAssign{#1/}{\the\pageno}}
\def\ix#1#2#3{\pAssign{#2}{\the\pageno}
\immediate\write#1{\noexpand\idxitem{#3}{\noexpand\csname+#2\endcsname}}}
\def\rf#1{\Wanted{#1}\csname+#1\endcsname\relax\rPut {#1}}
\def\rfp#1{\Wanted{#1}\csname+#1/\endcsname\relax\rPut{#1}}

%
\verno =2
\expandafter \def \csname +pone\endcsname {1}
\expandafter \def \csname +pr\endcsname {1.1}
\expandafter \def \csname +bc\endcsname {1.2}
\expandafter \def \csname +sla\endcsname {1.3}
\expandafter \def \csname +slb\endcsname {1.4}
\expandafter \def \csname +ec\endcsname {1.5}
\expandafter \def \csname +ptwo\endcsname {2}
\expandafter \def \csname +cohull\endcsname {2.1}
\expandafter \def \csname +sbg\endcsname {2.1}
\expandafter \def \csname +rem\endcsname {2.2}
\expandafter \def \csname +luno\endcsname {2.3}
\expandafter \def \csname +ldue\endcsname {2.4}
\expandafter \def \csname +pthree\endcsname {3}
\expandafter \def \csname +tolech\endcsname {3.1}
\expandafter \def \csname +setH\endcsname {3.1}
\expandafter \def \csname +aaa\endcsname {3.2}
\expandafter \def \csname +enea\endcsname {3.3}
\expandafter \def \csname +ene\endcsname {3.4}
\expandafter \def \csname +char\endcsname {3.2}
\expandafter \def \csname +pra\endcsname {3.5}
\expandafter \def \csname +prb\endcsname {3.6}
\expandafter \def \csname +prc\endcsname {3.7}
\expandafter \def \csname +prd\endcsname {3.8}
\expandafter \def \csname +pre\endcsname {3.9}
\expandafter \def \csname +prf\endcsname {3.10}
\expandafter \def \csname +prg\endcsname {3.11}
\expandafter \def \csname +ltre\endcsname {3.3}
\expandafter \def \csname +lqua\endcsname {3.4}
\expandafter \def \csname +coe\endcsname {3.12}
\expandafter \def \csname +esti\endcsname {3.13}
\expandafter \def \csname +res\endcsname {3.5}
\expandafter \def \csname +ref\endcsname {4}
\expandafter \def \csname +AAB\endcsname {1}
\expandafter \def \csname +AC\endcsname {2}
\expandafter \def \csname +BD\endcsname {3}
\expandafter \def \csname +BM\endcsname {4}
\expandafter \def \csname +CC\endcsname {5}
\expandafter \def \csname +CTZ\endcsname {6}
\expandafter \def \csname +ces\endcsname {7}
\expandafter \def \csname +Cl\endcsname {8}
\expandafter \def \csname +Cla\endcsname {9}
\expandafter \def \csname +Dei\endcsname {10}
\expandafter \def \csname +ET\endcsname {11}
\expandafter \def \csname +Mar\endcsname {12}
\expandafter \def \csname +Ol\endcsname {13}
\expandafter \def \csname +Ray\endcsname {14}

\immediate\openout1=\jobname.aux
\immediate\write1{\noexpand\verno=\the\verno}
\ifindexed
\immediate\openout2=\jobname.idx
\immediate\openout3=\jobname.sym
\fi
\def\prebye{\ifwanted
\message{Warning: Undefined references! Rerunning could help}\fi}
\headline={\ifnum\pageno>0\ifodd\pageno\rightheadline
\else\leftheadline\fi\fi}
\def\rightheadline{\hfil{\eightpoint\titolo}
\hfil\tenrm\folio}
\def\leftheadline{\tenrm\folio\hfil{\eightpoint\autori}
\hfil} \topskip=25pt
%
%
%

\def\R{{I\!\!R}}

\def \AC{\cal AC}

\def\F{{\cal F}}

%

%

%
%

\def\F{{\cal F}}

\def\AC{{\cal AC}}
\def\fs{f^{**}}

%
%
%
\def\titlea{An existence result for non-coercive non-convex}
\def\titleb{problems in the Calculus of Variations}
\def\author{GRAZIANO CRASTA}
\def\keywords{Calculus of Variations, existence theory,
non-coercive problems, non-convex problems, convex analysis.}
\def\abstract{
We consider the minimization problem $({\cal P})$:
$${\rm min}\ \int_0^T [a(t)\cdot u(t)+f(u'(t))]\,dt,$$
in the class $\AC$ of the absolutely continuous functions from $[0,T]$
into $\R^m$ satisfying the boundary conditions:
$$u(0)=u_0,\qquad u(T)=u_1.$$
We consider the class $\F$ of all lower semicontinuous functions
$f:\R^m\rightarrow\R$,
such that the convexification $\fs$ satisfies the growth condition:
$$\fs(x_n)-x_n\cdot\nabla\fs(x_n)\rightarrow -\infty,
\quad {\rm for}\ n\rightarrow+\infty,$$
for every sequence $(x_n)\subset\R^m$
of points of differentiability of $\fs$
such that
$|x_n|\rightarrow +\infty$ as $n\rightarrow+\infty$.
We prove that, for every $f\in\F$, the minimization problem
$({\cal P})$ has a solution.
}
\null
\pageno=0
\tenpoint
\baselineskip=15pt
{  
\leftskip=1.5truecm
\rightskip=2.5truecm
\vskip 2 truecm
\centerline{\medbf \titlea}
\vskip 4 truemm
\centerline{\medbf \titleb}
\vskip 10 truemm
\iftitlepage
\centerline{\author}
\vskip 3truemm
\centerline{\sl S.I.S.S.A, Via Beirut 2--4, 34014 Trieste (Italy)}
\centerline{\sl and}
\centerline{\sl Dip.~di Matematica Pura ed Applicata,
Via Campi 213/b, Modena (Italy)}
\vskip 10 truemm
\fi
%
\iftitlepage
 {\ninepoint 
   \noindent{\bf Abstract.\ }
   \abstract
 } 
 \vskip 1truecm
 \noindent{\bf Key Words.}\
 \keywords
 \vskip 1 truecm
 \centerline{Ref. S.I.S.S.A. 105/94/M, July 1994.}
\fi
}  
\iftitlepage
\vfill
\eject
\fi
%
%
%
\input olech1.tex
\vfill\eject
%
%
%
\newpar{ref}{References}

\bibart{AAB}{Ambrosio, L., Ascenzi, O.~and Buttazzo, G.}
{Lipschitz regularity for minimizers of integral functionals
with highly discontinuous integrands}
{J.~Math.~Anal.~Appl.}{142}{1989}{301}{316}

\biblib{AC}{Aubin, J.P.~and Cellina, A.}
{Differential Inclusions}
{Springer--Verlag}{Berlin}{1984}

\bibart{BD}{Botteron, B.~and Dacorogna, B.}
{Existence and non-existence results for non-coercive
variational problems and applications in ecology}
{J.~Diff.~Equations}{85}{1990}{214}{235}

\bibart{BM}{Botteron, B.~and Marcellini, P.}
{A general approach to the existence of minimizers of one
dimensional non-coercive integrals of the calculus of variations}
{Ann.~Inst.~H. Poin\-ca\-r\'e}{8}{1991}{197}{223}

\bibart{CC}{Cellina, A.~and Colombo, G.}
{On a classical problem of the calculus of variations
without convexity assumptions}
{Ann.~Inst.~H.~Poincar\'e}{7}{1990}{97}{106}

\bibprep{CTZ}{Cellina, A., Treu, G.~and Zagatti, S.}
{On the minimum problem for a class of non-coercive functionals}
{Preprint SISSA}{1994}

\biblib{ces}{Cesari, L.}
{Optimization -- theory and applications}
{Springer--Verlag}{New York}{1983}

\biblib{Cl}{Clarke, F.~H.}
{Optimization and Nonsmooth Analysis}
{Wiley Interscience}{New York}{1983}

\bibart{Cla}{Clarke, F.~H.}
{An indirect method in the calculus of variations}
{Trans.~Am.~Math. Soc.}{336}{1993}{655}{673}

\biblib{Dei}{Deimling, K.}
{Multivalued Differential Equations}
{De Gruyter}{Berlin}{1992}

\biblib{ET}{Ekeland, I.~and Temam, R.}
{Convex Analysis and Variational Problems}
{North--Holland}{Amsterdam}{1977}

\bibart{Mar}{Marcellini, P.}
{Alcune osservazioni sull'esistenza del minimo di integrali
del calcolo delle variazioni senza ipotesi di convessit\`a}
{Rendiconti di Matematica}{13}{1980}{271}{281}

\biblib{Ol}{Olech, C.}
{The Lyapunov theorem: its extensions and applications}
{in ``Methods of Non-convex Analysis'', Cellina Ed.}
{Springer--Verlag}{1990}

\bibart{Ray}{Raymond, J.P.}
{Existence theorems in optimal control problems without
convexity assumptions}
{J.~Optim.~Theory Appl.}{67}{1990}{109}{132}

\prebye
\bye